\documentclass[a4paper]{amsart}
\usepackage{amssymb}
\usepackage{amsmath,amsthm}
\usepackage{enumerate,paralist}
\usepackage{amstext}
\usepackage{psfrag}
\usepackage{dsfont}
\usepackage[dvips]{epsfig}
\usepackage[english]{babel}
\usepackage{url}
\usepackage{cite}
\usepackage{dcolumn}
\usepackage{bm}
\usepackage{float}
\usepackage{hyperref}
\usepackage{color}
\usepackage{epstopdf}
\usepackage{cleveref}
\usepackage{cite}
\usepackage[svgnames]{xcolor}
\hypersetup{hidelinks,colorlinks=true,allcolors=DarkBlue}
\usepackage{graphicx}
\theoremstyle{plain}
\newtheorem{theorem}{Theorem}

\newtheorem{lemma}{Lemma}
\theoremstyle{definition}

\newtheorem{problem}{Problem}

\def\sign{{\rm sign}}

\def\diag{{\rm diag}}

\numberwithin{equation}{section}

\begin{document}
\frenchspacing

\title[Feedback control for damping a system of linear oscillators]
{Feedback control for damping \\ a system of linear oscillators}

\author{Alexander Ovseevich}
\address
{
Institute for Problems in Mechanics, Russian Academy of Sciences \\
101/1 Vernadsky av., Moscow 119526, Russia.
} 
\email{ovseev@ipmnet.ru}

\author{Aleksey Fedorov}
\address
{
Institute for Problems in Mechanics, Russian Academy of Sciences \\ 101/1 Vernadsky av., Moscow 119526, Russia. \& 
Bauman Moscow State Technical University, 5 2nd Baumanskaya str., Moscow 105005, Russia.
} 
\email{akfedorov@student.bmstu.ru}

\maketitle

\begin{abstract}
The paper is devoted to a design of a common bounded feedback control steering a system of an arbitrary number of linear oscillators to the equilibrium. 
At high energies, the control is based on the asymptotic theory of reachable sets of linear systems. 
With decreasing of the energy, a similar control with a reduced upper bound is used. 
On the final stage, the control is constructed by using the method of common Lyapunov functions. 
Special attention is paid to the cases of one and two oscillators.


\end{abstract}

\section{Introduction}

In control theory, a class of problems related to the steering of a system from an arbitrary initial state to a given manifold in minimum time is well studied. 
Strictly speaking, this problem consists of two. The first of them is to bring the system to a given manifold. 
The second problem requires minimization of the motion time. 
There are two main approaches to design of the corresponding control. 
First, a feedback control (synthesis), where the control $u(x)$ is a function of the current state of the system. 
Second, a feedforward control $u(t)$, which is a function of the elapsed time. 
One of the classical results of control theory is the analytic design of the feedback  control for a single linear oscillator based on the Pontryagin maximum principle \cite{pont}.

\subsection{Problem statement}

In the present paper, a next in complexity problem is considered:
\begin{problem}\label{problem}
Design a feedback control for steering  a system of $N$ linear oscillators with eigenfrequencies $\omega_i$
\begin{eqnarray}\label{syst1}
	&\dot{x}={A}x+{B}u,\quad x=(x_1,y_1,\dots,x_N,y_N)^*\in{\mathbb{R}}^{2N},\quad u\in {\mathbb{R}},\quad |u|\leq1,\\[1em]
\label{syst2} &{A} = \left( {\begin{array}{*{20}c}
	0	&	1	&	&	&	\\
	{-\omega_1^{2}}	&	0	&	&	&	\\
	&	&	\ddots	&	&	\\
	&	&	&	0	& 	1	\\
	&	&	&	{-\omega_N^{2}}	&	0	\\
\end{array} } \right), \qquad
	{B} = \left(\begin{array}{c}
	0	\\
	1	\\
	\vdots	\\
	0	\\
	1	\\
\end{array}\right),
\end{eqnarray}
to the equilibrium in minimum time.
\end{problem}
Here and in what follows $*$ stands for the transposition.

It looks like that an analytic construction of the optimal feedback control by using maximum principle methods in Problem \ref{problem} is impossible, 
and even the search for a numerical solution is a hard problem.

\begin{figure}
\begin{minipage}[h]{0.45\linewidth}
	\center{\resizebox{1.0\columnwidth}{!}{\includegraphics{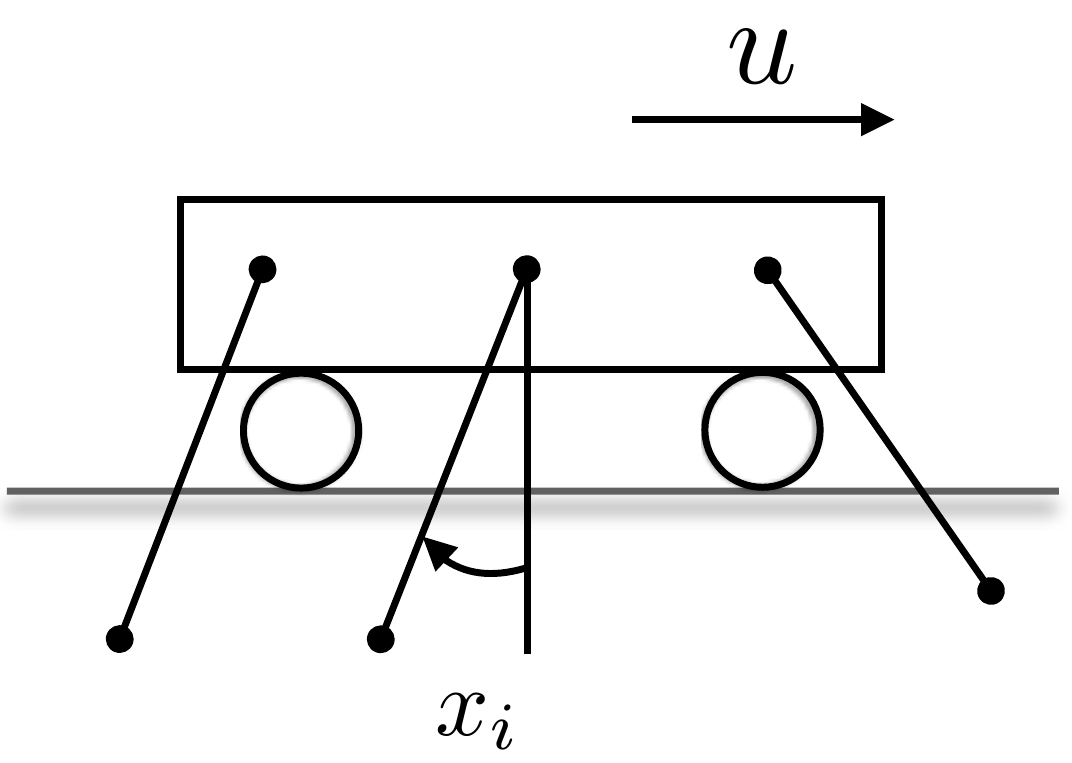}}} \mbox{a)}
\end{minipage}
\hfill
\begin{minipage}[h]{0.45\linewidth}
	\center{\resizebox{1.0\columnwidth}{!}{\includegraphics{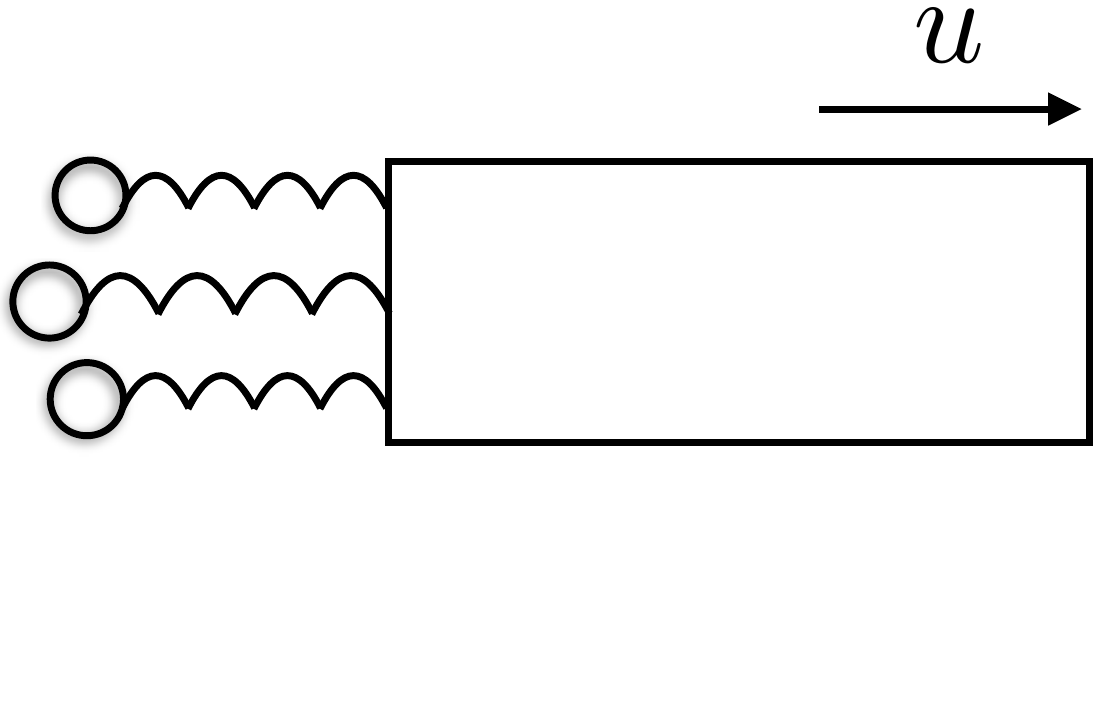}}} \mbox{b)}
\end{minipage}
\hfill 
\caption
{ 
	A mechanical model of the system: 
	a) system of $N$ pendulums with eigenfrequencies $\omega_i$ attached to a cart; 
	b) masses, attached via springs to a body, which moves with acceleration. 
} 
\label{fig:system}
\end{figure}

A mechanical model for system (\ref{syst1})--(\ref{syst2}) is a system of $N$ pendulums with eigenfrequencies $\omega_i$ attached to a cart moving with bounded acceleration $u$, 
where the vertical deviations $x_i$ of pendulums are components of the state vector (Fig. \ref{fig:system}a). 
In the other interpretation, components of the vector are the displacements $x_i$ of masses, attached to a body moving with acceleration $u$ (Fig. \ref{fig:system}b).

\subsection{Proposed approach}

A non-optimal feedback control, which steers the system to the equilibrium, will be used. 
The control is asymptotically optimal: 
a ratio the time motion to zero by using this control to the minimum one is close to one if the initial energy of the system
\begin{equation}
	E=1/2\sum_{i=1}^{N}{(\dot{x_i}^2+\omega_i^2 x_i^2)}
\end{equation}
is sufficiently large. 
Proof of the asymptotic optimality of the suggested control is beyond our paper and it will be presented elsewhere.

In our paper, an assumption of absence of resonance, {\it i.e.}, non-trivial relations between eigenfrequencies of the form,
\begin{equation}\label{reson}
	\sum_{i=1}^Nm_i\omega_i=0,\, \mbox{ where } 0\neq m=(m_1,\dots,m_N)\in {\mathbb{Z}}^N
\end{equation}
is used. 
In particular, this implies that for system (\ref{syst1})--(\ref{syst2}) the Kalman controllability condition holds \cite{kalman}. 
The suggested control works also in the resonant case, however its quasi-optimal properties are lost.

\begin{figure}
\begin{minipage}[h]{0.3\linewidth}
	\center{\resizebox{1.0\columnwidth}{!}{\includegraphics{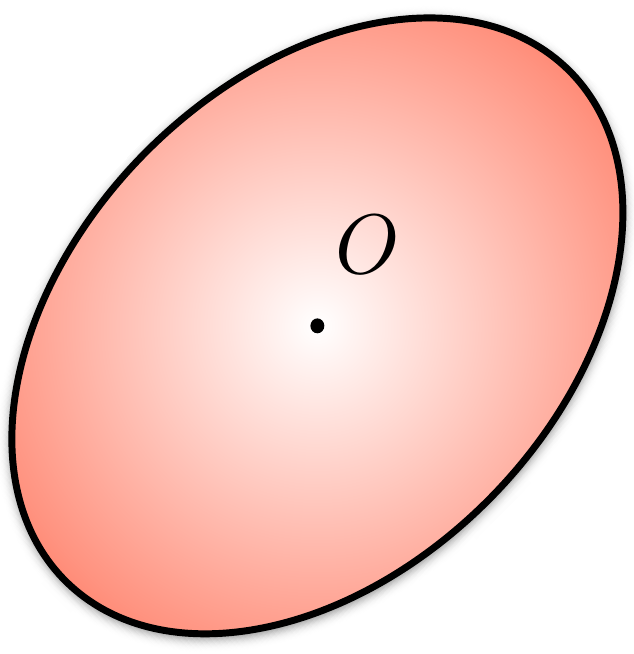}}} \mbox{a)}
\end{minipage}
\hfill
\begin{minipage}[h]{0.3\linewidth}
	\center{\resizebox{1.0\columnwidth}{!}{\includegraphics{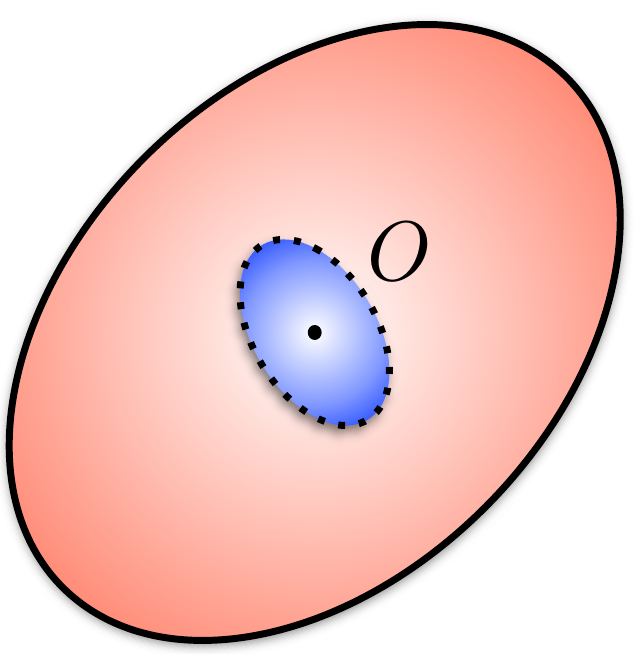}}} \mbox{b)}
\end{minipage}
\hfill
\begin{minipage}[h]{0.3\linewidth}
	\center{\resizebox{1.0\columnwidth}{!}{\includegraphics{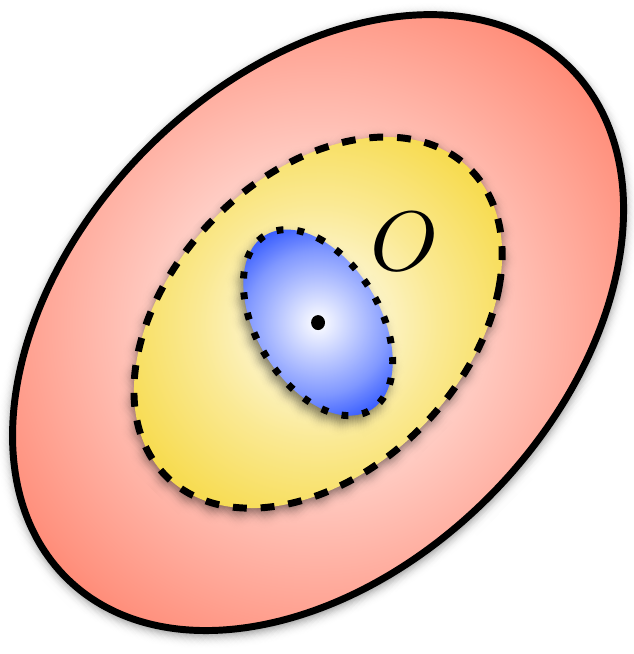}}} \mbox{c)}
\end{minipage}
\caption
{
	Structure of the suggested control: 
	the region bounded by the solid line is the region in which the control based on the asymptotical theory of reachable sets is used; 
	the region bounded by the dashed line is the region in which analogical control with reduced upper bound is used; 
	the region bounded by the dotted line is the region in which the local feedback is applied. 
	The point $O$ is the equilibrium position.
} 
\label{fig:2}
\end{figure}

To solve the  problem, three strategies are used in sequence.
\begin{enumerate}
	\item {\it High energy zone}.
	At high energies, the normal vector to the {\it approximate} reachable set, close to the exact one in the long run, serves as a momentum $p(x)$ (Fig. \ref{fig:2}a) \cite{ovseev6,ovseev,ovseev2}.
	\item {\it Intermediate zone}.
	The obtained control is applicable also at low energies, but its quasi-optimal properties in this case are lost.
	First, this control affects  the system like a dry friction, so that in some states it prohibits any motion.
	Second, the motion may occur in a small neighborhood of a limit set (attractor), which does not contain the equilibrium position.
 	Application of the control with reduced upper bound (second stage of the control) allows to postpone this undesirable pulling into the attractor.
	\item {\it Low energy zone}.
	At the final third stage, an approach to the design of a local control in the feedback form based on the common Lyapunov functions is used.
	This approach works in a small neighborhood of the equilibrium position (Fig. \ref{fig:2}b).
	To reach this small neighborhood, it should necessarily  contain the basin of attractor  of the preceding control.
	The reduction of the basin achieved at the previous stage makes this objective possible (Fig. \ref{fig:2}c).
\end{enumerate}

It should be noted that there exist other methods of control design in the feedback form, 
{\it e.g.}, based on the Kalman approach to the feedforward control for linear systems \cite{chern,ovseev3}. 
The cited works also include estimates for the motion time under the control, 
which are comparable with the optimal time in the sense that the ratio the motion time under the  the proposed control to the minimum one is bounded.

\subsection{Goal of the work}

A possibility of a design containing three steps for the system of linear oscillators has been demonstrated for the first time in Ref. \cite{ovseev4}. 
A brief summary of results on the asymptotic optimality of the suggested control, partially based on a new lemma on the observable time-invariant linear systems, 
and on the study of the character of singular motion, can be found in Ref. \cite{ovseev5}.

In the present paper, we emphasize the computational and ``algorithmic'' aspects of the suggested control. 
The paper is organized as follows. 
Section \ref{high} is devoted to a design, computation, and properties of the control used in high energy and intermediate zones. 
We give the explicit expressions for the support function of the reachable set of system (\ref{syst1})--(\ref{syst2}). 
In the general case, the support function is expressed via the Gelfand hypergeometric function. 
In the case of two oscillators it is expressed via the elliptic integrals. 
Methods for solution of the master equation determining the control are indicated. 
In the case of two oscillators, the method is based on the solution of a transcendental equation in elliptic functions. In the case of an arbitrary $N$, we state an optimization problem such
that its solution is equivalent to the solution of the master equation. 
For the case of two oscillators $N=2$ the solution of the optimization problem is equivalent to the solution of a transcendent equation in elliptic functions.
Section \ref{terminal} is devoted to construction of the control in a region close to the terminal point. 
In particular, common Lyapunov functions for $N=2$ and $N=1$ are calculated, it is shown also that for the arbitrary $N$ the common
Lyapunov function is an even integer-valued  matrix. In section \ref{comb}, main matching relations for the control at different stages are presented. 
In Section \ref{toy}, the suggested strategy is illustrated for the case of a single linear oscillator.

Proofs of all the statements can be found in \cite{fedorov}.

\section{Control at high energies}\label{high}

\subsection{Asymptotic theory of reachable sets}

Being simultaneously an instrument and an object of studies, reachable sets occupy a central place in control theory.
Therefore, their study is a fundamental problem \cite{chern2}. 
Recall that the reachable set $\mathcal{D}(T)$ is defined as the set of points reachable from zero at time $T$.

One of the main result of Ref. \cite{ovseev} as applied to the considered system (\ref{syst1})--(\ref{syst2}) is as follows: 
the reachable set $\mathcal{D}(T)$ has as $T\to\infty$ asymptotic form $T\Omega$, where $\Omega$ is a fixed convex body.
The result is stated more precisely in the following theorem.

\begin{theorem}\label{support0}
	Suppose that the momentum $p$ is written in the form $p=(\xi_i,\eta_i)$,
	where
	$\xi_i$ is the dual variable for $x_i$,
	$\eta_i$ is the dual variable for $y_i$, and
	$z_i=(\eta_i^2+{\omega_i^{-2}}{\xi_i^{2}})^{1/2}$.
	In the case of absence of resonance (nontrivial relations between eigenfrequencies of the form (\ref{reson})),
	the support function $H_T$ of the reachable set $\mathcal{D}(T)$ has as $T\to\infty$ the asymptotic form
	\begin{equation}\label{approxN0}
		{H}_T(p)=\frac{T}{(2\pi)^N}\int_{0}^{2\pi}\dots\int_{0}^{2\pi}\left|\sum_{i=1}^N z_{i}\cos\varphi_{i}\right|d\varphi_{1}\dots d\varphi_{N}+o(T)=T\mathfrak{H}(z)+o(T),
	\end{equation}
	and the support function of the compact set $\Omega$ is given by the main term $\mathfrak{H}(z)$.
\end{theorem}

The support function of any subset $M\in\mathbb{R}^{N}$ is defined as
\begin{equation}\label{SF}
	{H}_{M}(\xi)=\sup_{x\in M}\langle{\xi, x\rangle},
\end{equation}
where angle brackets stand for the scalar multiplication in $\mathbb{R}^{N}$.

\subsection{Approach to the control design}

Geometrically, 
the Pontryagin maximum principle says that the momentum (vector of adjoint variables) $\psi$ in a point $x$ is the inner normal to the reachable set $\mathcal{D}(T)$ 
($T(x)$ being the time of reaching $x$ from zero).

The idea of our approach to control design is to use the set $T\Omega$ as an approximation to the reachable set $\mathcal{D}(T)$, and the normals to $T\Omega$ as momenta. 
If the phase vector $x\in{V}=\mathbb{R}^{2N}$ lies at the boundary of the set $T\Omega$, then
\begin{equation}\label{approx3}
	{T}^{-1}x=\frac{\partial {H}_\Omega(p)}{\partial p}
\end{equation}
for a certain momentum $p=p(x)$. 
We note that the support function ${H}_{\Omega}$ is differentiable at $p\ne0$ and Eq. (\ref{approx3}) 
has a unique solution because of the smoothness of the boundary of the set $\Omega$ proved in Ref. \cite{ovseev2}. 
The feedback control is given by
\begin{equation}\label{approx4}
	{u(x)}=-\sign\langle{B,p(x)\rangle}.
\end{equation}

In what follows, we will also use the control in the form
\begin{equation}\label{controlU}
	u_U(x)=Uu(x),\quad |U|\leq1.
\end{equation}

\subsection{Computation of the control}

In the coordinates notation, Eq. (\ref{approx3}) has the form
\begin{equation}\label{approx_coord}
	T^{-1}(x_i,y_i)=z_i^{-1}\left(\frac{\partial {\mathfrak H}}{\partial z_i}\right)\left(\frac{\xi_i}{\omega_i^2},\eta_i\right), \quad i=1,\dots,N,
\end{equation}
where $z_i=(\eta_i^2+{\omega_i^{-2}}{\xi_i^{2}})^{1/2}$. 
In order to solve Eq. (\ref{approx_coord}) we  calculate first the point $\mathfrak z$ of a sphere $\mathcal{S}^{N-1}$ with positive-homogeneous coordinates $(z_1:\dots:z_N)$. 
Here the sphere $\mathcal{S}^{N-1}$ is regarded  as a set of directions of non-zero vectors in $\mathbb{R}^N$. 
To this end, 
define an ``energetic'' vector $e=(e_i)\in\mathbb{R}^N$, where  $e_i=(\omega_i^{2}x_i^{2}+y_i^{2})^{1/2}$ and obtain from Eq. (\ref{approx_coord}) that
\begin{equation}\label{z}
	T^{-1}e_i=\frac{\partial {\mathfrak H}}{\partial z_i},\quad i=1,\dots,N.
\end{equation}

The solution of Eq. (\ref{approx_coord}) defines an inversion of a map from one $2N$-dimensional manifold to another one,
while the solution of (\ref{z}) reduces to inversion of map of $(N-1)$-dimensional manifolds. 
Like Eq. (\ref{approx_coord}), the Eq. (\ref{z}) possesses a unique solution \cite{ovseev2}, which can found as follows.

In the case $N=1$, the support function has the form $\mathfrak{H}(z)=\frac2\pi|z|$, and the solution is ${\mathfrak z}=e$, 
while the control has the dry-friction form $u=-\sign y_1$. 
If $N=2$, then Eq. (\ref{z}) is one-dimensional: 
the function $\mathfrak{H}$ can be expressed through elliptic integrals, point $(e_1:e_2)\in\mathcal{S}^1$ defines problem (\ref{z}) completely, 
and the solution $\mathfrak z\in\mathcal{S}^1$ can be also found by solving a transcendental equation in elliptic functions.

\subsubsection{Formula for the support function}
The support function of the convex body $\Omega$ is given by the main term of the asymptotic expression in
(\ref{approxN0})
\begin{equation}\label{approx2N}
	{H}_\Omega(p)=\mathfrak{H}(z)=\int\left|\sum_{i=1}^N z_{i}\cos\varphi_{i}\right|d\varphi, \mbox{ where } z=(z_1,\dots,z_N)\in{\mathbb R}^N.
\end{equation}
For $N=1$, we obtain $\mathfrak{H}(z)=\frac2\pi|z|$, in the case $N=2$ the function
\begin{equation}
	\mathfrak{H}(z)=\int\left|z_{1}\cos\varphi_{1}+z_{2}\cos\varphi_{2}\right|d\varphi
\end{equation}
can be expressed via elliptic integrals (we will dwell on this subject below). In the general case,
\begin{equation}
	\mathfrak{H}(z)=\frac{1}{(2\pi)^N}\int\limits_{\{|t_i|\leq1\}}{\left|\sum_{i=1}^N z_{i}t_{i}\right|}{\prod_{i=1}^N(1-t_i^2)^{-1/2}}dt_1\dots t_N
\end{equation}
is an Euler-type integral that defines a generalized hypergeometric function in the sense of I.M. Gelfand

\subsubsection{Case of two oscillators}

Our previous arguments imply that the case $N=2$ is special in some sense, because in this case one can solve Eq. (\ref{approx3}) by {\it a simpler} method. 
Indeed, in this case
\begin{equation*}\label{ellptic1}
	\frac{\partial {\mathfrak H}}{\partial z_i}=\frac{1}{2\pi}\int_{0}^{2\pi} \cos\phi_i \sign(z_{1}\cos\varphi_{1}+z_{2}\cos\varphi_{2})d\varphi.
\end{equation*}

For instance, consider the index $i=1$ and make the inner integration wrt  $\varphi_{2}$. 
Taking into account the non-negativity of $z_2$, one needs to compute the integral
\begin{equation}\label{ellptic3}
	\frac{1}{2\pi} \int_0^{2\pi}\sign(-C+\cos\varphi_{2})d\varphi_{2}=\frac{2}{\pi}\arccos C-1, \mbox{ если }|C|\leq1,
\end{equation}
where $C=k\cos\phi_1,$ $k=-z_1/z_2$. It is obvious that by permutation of the indexes, we can assume that $|k| \leq 1$.
Note that the latter assumption breaks the symmetry between $z_1$ and $z_2$. At $|k|\leq1$, from (\ref{ellptic3}) we obtain that
\begin{equation}\label{ellptic4}
	\frac{\partial {\mathfrak H}}{\partial z_1}=\frac{1}{\pi^2}\int_0^{2\pi}\cos\varphi_1\arccos (k\cos\varphi_1)d\varphi_1,
\end{equation}
because $\int_0^{2\pi}\cos\varphi_1 d\varphi_1=0$. 
After integration by parts, integral (\ref{ellptic4}) can be written in the elliptic form
\begin{equation}\label{ellptic5}
	\int_0^{2\pi}\cos\varphi\arccos (k\cos\varphi)\varphi=\int_0^{2\pi}\frac{k\sin^2\varphi}{\sqrt{1-k^2\cos^2\varphi}}\,d\varphi
\end{equation}
and finally we get for the following formula for the derivative of the support function
\begin{equation}\label{ellptic6}
	\frac{\partial {\mathfrak H}}{\partial z_1}=\frac{1}{\pi^2}\int_0^{2\pi}\frac{k\sin^2\varphi}{\sqrt{1-k^2\cos^2\varphi}}\, d\varphi,
	\mbox{ where }k=-z_1/z_2,
\end{equation}
which holds true provided that $|k|\leq1$. 
For computation of $\frac{\partial {\mathfrak H}}{\partial z_2}$, we need to compute the inner integral
\begin{equation*}\label{ellptic31}
	\frac{1}{2\pi} \int_0^{2\pi}\cos\varphi_{2}\sign(-C+\cos\varphi_{2})d\varphi_{2}=\frac{2}{\pi}\sin\arccos C, \mbox{ if }|C|\leq1,
\end{equation*}
wherefrom we obtain
\begin{equation}\label{ellptic7}
	\frac{\partial {\mathfrak H}}{\partial z_2}=\frac{1}{\pi^2}\int_0^{2\pi}{\sqrt{1-k^2\cos^2\varphi}}\, d\varphi.
\end{equation}
Note that the asymmetry between integral formulas (\ref{ellptic6}) and (\ref{ellptic7}) is illusive: 
the change of variables $z_1\leftrightarrows z_2$ implies the change of parameters $k\leftrightarrows k^{-1}$. 
Under this change, the integrals
\begin{equation}
	I_1(k)=\int_0^{2\pi}\frac{k\sin^2\varphi}{\sqrt{1-k^2\cos^2\varphi}}\, d\varphi\mbox{ and }I_2(k)=\int_0^{2\pi}{\sqrt{1-k^2\cos^2\varphi}}\, d\varphi,
\end{equation}
regarded as meromorphic functions of  $k$, are transposed: $I_1(k^{-1})=I_2(k)$.

The key equation (\ref{z}) that defines control (\ref{approx4}), has the form of equation for $k=-z_1/z_2$:
\begin{equation}\label{ellptic_control}
	\frac{e_2}{e_1}=\frac{I_2}{I_1}{(k)}.
\end{equation}

For numerical experiments, it is useful to employ a representation for (\ref{ellptic4}) and (\ref{ellptic7}) via canonical Legendre elliptic integrals of the first and second kinds
\begin{equation}
	\mathcal{A}(k)=\sqrt{1-k^2}\mathcal{E}\left(\frac{k^2}{k^2-1}\right), \quad
	\mathcal{E}\left(y\right)=\int_0^{\pi/2}{\sqrt{1-y^2\sin\varphi}\,d\varphi},
\end{equation}
where $\mathcal{E}\left(y\right)$ is the complete Legendre integral of the second kind,
\begin{equation}
	\mathcal{B}(k)=-\frac{\sqrt{1-k^2}}{k}\left[\mathcal{E}\left(\frac{k^2}{k^2-1}\right)-\mathcal{K}\left(\frac{k^2}{k^2-1}\right)\right], \quad
	\mathcal{K}\left(y\right){=}{\int}_0^{\pi/2}\frac{d\varphi}{\sqrt{1-y^2\sin\varphi}},
\end{equation}
where $\mathcal{K}\left(y\right)$ is the complete Legendre integral of the first kind.

Then expression (\ref{ellptic_control}) can be rewritten in the form
\begin{equation}\label{ellptic_control2}
	F(k)={\frac{\partial{\mathfrak H}}{\partial z_2}}\bigg\slash{\frac{\partial{\mathfrak H}}{\partial z_1}}=\frac{\mathcal{B}(k)}{\mathcal{A}(k)}.
\end{equation}
The graph of the function $F(k)$ is presented in Fig. \ref{fig:2}. 
Note that the support function has the form
\begin{equation}\label{ellptic_integral}
	\mathfrak{H}(z_1,z_2)=\frac{1}{\pi^2}\int_0^{2\pi}\frac{(z_2^2-z_1^2)d\varphi}{\sqrt{z_2^2-z_1^2\cos^2\varphi}}\quad \mbox{ для } |z_1|\leq|z_2|.
\end{equation}

Note that the function ${\partial{\mathfrak H}}/{\partial z_i}$ can be regarded as a function  $g_i=g_i(\mathfrak z)$ on $\mathcal{S}^{N-1}$. 
Assuming that $g_i$ are known, we obtain the final formula for the momentum
\begin{equation}\label{impulse}
	({\xi_i},\eta_i)=\frac{z_i}{Tg_i}({\omega_i^2}x_i,y_i),\, i=1,\dots,N.
\end{equation}
If we know the point $\mathfrak {z}=(z_1:\dots:z_N)\in\mathcal{S}^{N-1}$, 
then due to the non-negativity of the coordinates $z_i$ the direction of the momentum $p(x)$ is defined by formula (\ref{impulse}) uniquely: 
the unknown positive pre-factor $T^{-1}$ does not play a role. 
Control (\ref{approx4}) depends on the direction of the momentum only, and therefore it can be efficiently found. It has the form
\begin{equation}\label{control}
	u(x)=-\sign\left(\sum_{i=1}^N{g_i}^{-1}{z_iy_i}\right).
\end{equation}

\begin{figure}
\begin{center}
\begin{minipage}[h]{0.75\linewidth}
\center{\resizebox{1.0\columnwidth}{!}{\includegraphics{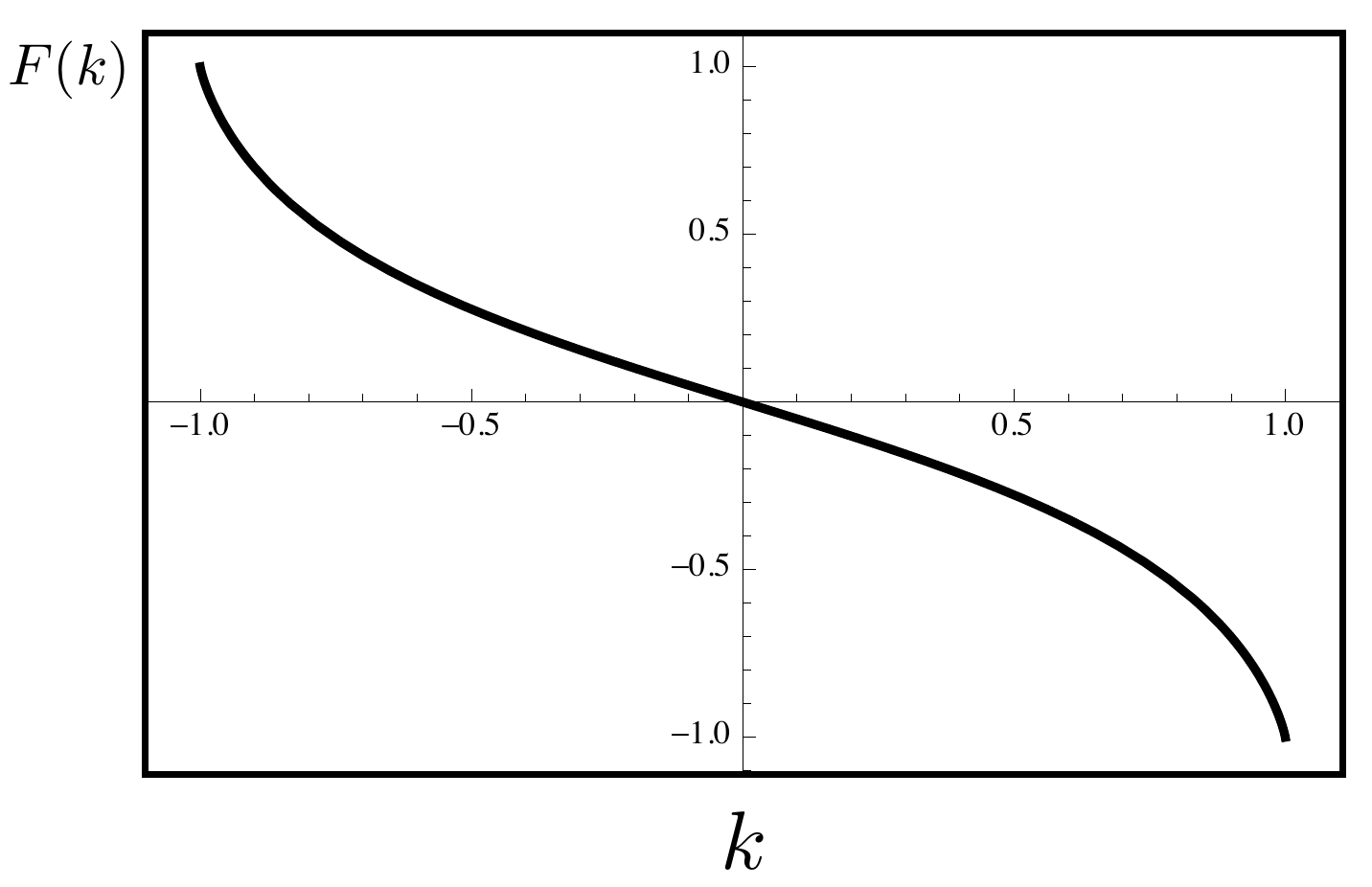}}}
\end{minipage}
\end{center}
\caption
{
	Plot of the function $F(k)$ governed by Eq. (\ref{ellptic_control2}).
} 
\label{fig:3}
\end{figure}

\subsubsection{Case of an arbitrary number of oscillators}

For an arbitrary number of linear oscillators, thanks to the Kuhn--Tucker theorem the search for a solution of Eq. (\ref{z}) is equivalent to solving the following optimization problem
\begin{equation}\label{optim}
    \langle{e,z}\rangle\to\max \mbox{, provided that }\mathfrak H(z)\leq1.
\end{equation}
It is obvious that the restriction $\mathfrak H(z)\leq1$ is equivalent to the restriction $\mathfrak H(z)=1$. The
hypersurface $\{\mathfrak H(z)=1\}$ is strictly convex due to the relation
\begin{equation}\label{hessian}
     \left\langle{\frac{\partial^2 {\mathfrak H}}{\partial z^2}(z)\xi,\xi}\right\rangle=\int\limits_V\left(\sum_{i=1}^N\xi_{i}\cos\varphi_{i}\right)^2d\sigma(\varphi),
\end{equation}
where the integration is over $V=\left\{z\in\mathbb{R}^N:f(z)=0\right\}$ and
\begin{equation}
    f(z)=\sum_{i=1}^N z_{i}\cos\varphi_{i},\quad d\sigma(\varphi)=\frac{d\varphi_1\wedge\dots\wedge d\varphi_N}{(2\pi)^N df}
\end{equation}
is the canonical volume element on $V$.

It follows from relation (\ref{hessian})  that if the vectors $\xi$ and $z$ are not collinear, then
\begin{equation}
	\left\langle{\frac{\partial^2 {\mathfrak H}}{\partial z^2}(z)\xi,\xi}\right\rangle>0.
\end{equation}
However if the vector $\xi$ is tangent to the hypersurface $\{\mathfrak H(z)=1\}$ in the point $z$, then these two vectors cannot be collinear. 
Indeed, the tangency condition $\left\langle{{\partial\mathfrak H}/{\partial z},z}\right\rangle=0$ contradicts the Euler relation 
$\left\langle{{\partial\mathfrak H}/{\partial z},z}\right\rangle={\mathfrak H}(z)>0$ for the homogeneous function ${\mathfrak H}$. 
Strict convexity of the hypersurface $\{\mathfrak H(z)=1\}$ means uniqueness of the solution of the optimization problem (\ref{optim}). 
From this strict convexity it follows that the function $f={\mathfrak H}^2$ is strictly convex as well.

In the same time, the optimization problem (\ref{optim}) is equivalent to
\begin{equation}\label{optim2}
	\langle{e,z}\rangle\to\max \mbox{, provided that } f(z)\leq1.
\end{equation}
if $z_1\neq z_2$ are solutions of (\ref{optim2}), then $\langle{e,z_1}\rangle=\langle{e,z_2}\rangle$ and $f(z_i)=1$.  
This means, however, that 
\begin{equation}
	\left\langle{e,\frac{z_1+z_2}{2}}\right\rangle=\langle{e,z_1}\rangle 
	\mbox{ and } 
	{f}\left(\frac{z_1+z_2}{2}\right)<1, 
\end{equation}
which contradicts  the optimality of vectors $z_i$.

Thus, the equivalent optimization problems (\ref{optim}) and (\ref{optim2}) can be solved by using the well-known efficient methods, {\it e.g.}, using the Matlab Optimization Toolbox.

\subsection{Asymptotic optimality of the control}

Define a polar-like coordinate system, which is well suited for description of the motion under the control $u$. 
Write the phase vector $x$ in the form $x=\rho\phi$, where $\rho>0$ and $\omega=\partial\Omega$. 
In the terms of Eq. (\ref{approx3}), $\rho=T$, and $\phi={\partial H_\Omega}/{\partial p}$. 
In these coordinates, the equation of motion has the form
\begin{equation}\label{T}
	\dot{\rho}=-\left|{\left\langle\frac{\partial{\rho}}{\partial x},B\right\rangle}\right|, 
	\quad
	\dot{\phi}=A\phi+\frac{1}{\rho}\left(Bu+\phi \left|\left\langle\frac{\partial {\rho}}{\partial x},B\right\rangle\right|\right).
\end{equation}
If $N=1$, we get the proper polar coordinate system in the plane.

For the function $\rho=\rho(x)$, an eikonal-type equation holds
\begin{equation}\label{Euler}
	H_\Omega\left(p\right)=1,\quad p=\frac{\partial \rho}{\partial x}.
\end{equation}
It is dual to the equation $\rho({\partial {H_\Omega}}/{\partial p})=1$ of the surface $\omega$. 
Eq. (\ref{Euler}) can be used for averaging  the right-hand side of the first identity in (\ref{T}) with respect to time, 
and the proof of the following statement about the asymptotic optimality of control (\ref{approx4}) is based on the equation:
\begin{theorem}\label{main_approx}
	Consider the evolution of the value $\rho$ under control (\ref{approx4}). 
	Let $M=\min\{\rho(0),\rho(T),T\}$. 
	Then as $M\to+\infty$ we have
	\begin{equation}\label{approx_T}
		{(\rho(0)-\rho(T))}/{T}=1+o(1).
	\end{equation}
	Under any other admissible control,
	\begin{equation}\label{approx_T2}
		{(\rho(0)-\rho(T))}/{T}\leq 1+o(1).
	\end{equation}
\end{theorem}

The proof of this statement is given in Ref. \cite{fedorov}.

\subsection{Comparison with the maximum principle}

Solution of the linear minimum-time problem reduces completely to the boundary value problem for the Pontryagin maximum principle corresponding to the Hamiltonian
\begin{equation}
	h(x,\psi)=\langle{Ax,\psi}\rangle+|B^*\psi|-1=\max\{\langle{Ax,\psi}\rangle+\langle{Bu,\psi}\rangle-1\},
\end{equation}
where $|\cdot|$ is the Euclidean norm, and the maximum is taken over the interval $\{u\in\mathbb{R}:|u|\leq1\}$. 
The problem has the form
\begin{equation}\label{max}
\begin{array}{l}
	\dot{x}={A}x+{B}u,  \, \dot{\psi}=-{A^*}\psi,\\[.5em]
	u=\sign\langle{B,\psi}\rangle,\, x(0)=x_0,\,x(T)=0, \, h(x,\psi)=0.
\end{array}
\end{equation}

One can approach the issue of optimality of control (\ref{approx4}) by comparing the differential equations of the motion under the control with equations of the Pontryagin maximum principle. 
To this end, one need to understand how the momentum $p(x)$ in the Eq. (\ref{approx4}) changes with time. 
This description is given by the following equation
\begin{equation}\label{attractor_syst22p1}
	\dot{p}=-A^*p+\widetilde Bu, \mbox{ where } \widetilde{B}=\frac{\partial^2\rho}{\partial x^2}B.
\end{equation}
Note that if the second term $\widetilde{B}u$ is absent from the latter equation, 
then the equation for $\psi=-p$ would coincide with the equation of the maximum principle for the conjugate variables. 
However, the matrix $\frac{\partial^2\rho}{\partial x^2}$ is a homogeneous function of $x$ of degree $-1$, 
and, then the mentioned second term has the order of magnitude $O(1/{|x|})$ at large $x$ so that it is small. 
Also note that the condition of maximum $u=\sign\langle{B,\psi}\rangle=-\sign\langle{B,p}\rangle$ holds for control (\ref{approx4}). 
The motion under control (\ref{approx4}) is governed by the Hamiltonian 
$\mathcal{H}=\langle{Ax,\psi}\rangle+|\langle{B,\psi}\rangle|-\left|\left\langle{B,\partial\rho/\partial x}\right\rangle\right|$, which is in some sense close to the Pontryagin Hamiltonian $h(x,\psi)$. 
The difference between the Hamiltonians $\mathcal{H}$ and $h$ is $1-\left|\left\langle{B,\partial\rho/ \partial x}\right\rangle\right|$. 
The average value $\left|\left\langle{B,\partial\rho/\partial x}\right\rangle\right|$ is close to one for sufficiently large $x$, which follows from Theorem \ref{main_approx}.

Thus, for the vector $(x,\psi)$, where $\psi=-\partial\rho/\partial x$ the equation of the maximum principle holds ``on average'' with a small error at large $x$.

\section{Feedback near the terminal point}\label{terminal}

The idea of a design of local feedback control, which is a base of the third stage of the suggested control, 
goes back to \cite{anan} and uses a preliminary reduction of system (\ref{syst1})--(\ref{syst2}) to a canonical form by using transformations
\begin{equation}\label{transformations}
	A\mapsto A+BC,\quad  u\mapsto u-Cx ,\quad A\mapsto D^{-1}AD,\quad B\mapsto D^{-1}B,
\end{equation}
that correspond to adding a linear feedback control, and to coordinate changes. 
Formulate the result as a Lemma.

\begin{lemma}\label{canonical}
	By using transformations (\ref{transformations}), system (\ref{syst1})--(\ref{syst2}) reduces to the form:
	\begin{equation}\label{AB0}
		\dot {\mathfrak x}={\mathfrak A}{\mathfrak x}+{\mathfrak B}{\mathfrak u},
	\end{equation}
	\begin{equation}\label{AB}
	\begin{array}{c}
		\mathfrak{A} = \left( {\begin{array}{ccccc}
		0	&	&	&	\\
		-1	&	0	&	&	\\
		&	-2	&	0	&	\\
		&	&	\ddots	&	\ddots	\\
		&	&	&	-2N+1	&	0	\\
	\end{array} } \right), \quad
		\mathfrak{B} = \left( \begin{gathered}
		1 \hfill \\
		0 \hfill \\
		0 \hfill \\
		\vdots \hfill \\
		0 \hfill \\
		\end{gathered}  \right).
		\\
	\end{array}
	\end{equation}
Wherein the matrix of the linear feedback has the form
\begin{equation}\label{C}
	C=(c_1\,0\,c_2\,0\,\dots\, c_N\,0),\,c_k=(-1)^{N+1}\omega_k^{2N}\prod_{i\neq k}(\omega_i^{2}-\omega_k^{2})^{-1}.
\end{equation}
The gauge matrix $D$ has the following form. Define $2\times2$ matrices
\begin{equation}\label{d}
	d_{ij}=(-1)^{j-1}\lambda_i^{j-1}\left(%
 	\begin{array}{cc}
		0 &  -\frac{1}{(2j-1)!} \\
		\frac{1}{(2(j-1))!} & 0
	\end{array}\right), \mbox{ where }
	\lambda_k=\sum_{i\neq k}\omega_i^2.
\end{equation}
Then, D is the $N\times N$ matrix $(d_{ij})$ of $2\times2$ blocks $d_{ij}$.
\end{lemma}

As an existence theorem of a canonical form without explicit formulas for matrices $C$ and $D$, Lemma \ref{canonical} is a particular case of the Brunovsky theorem \cite{brun}. 
By following \cite{anan}, we introduce a matrix function of time related to system (\ref{AB}):
\begin{equation}\label{delta}
	\delta(\mathfrak T)=\diag({\mathfrak T}^1,{\mathfrak T}^{2}, \dots, {\mathfrak T}^{2N})^{-1}.
\end{equation}
In what follows the parameter ${\mathfrak T}$ will become a function ${\mathfrak T}={\mathfrak T}(\mathfrak{x})$ of the phase vector. 
Define the matrices in accordance with \cite{anan,korob}
\begin{equation}\label{M}
\begin{array}{l}
	\mathfrak{q}=(\mathfrak{q}_{ij}),\, \mathfrak{q}_{ij}=\int_0^1 x^{i+j-2}(1-x)dx=[(i+j)(i+j-1)]^{-1}, \\[1em]
	\mathfrak{Q}=\mathfrak{q}^{-1},\quad
	\mathfrak{C}=-\frac{1}{2}\mathfrak{B}^{*}\mathfrak{Q}, \quad \mathfrak{M}=\diag(1, 2, \dots, 2N). \\
\end{array}
\end{equation}
Define the feedback control by the equation
\begin{equation}\label{u}
	{\mathfrak u}(\mathfrak{x})=\mathfrak{C}\delta({\mathfrak T}( \mathfrak{x}))\mathfrak{x},
\end{equation}
where the function ${\mathfrak T}={\mathfrak T}(\mathfrak{x})$ is defined implicitly by the following relation:
\begin{equation}\label{condu}
	\langle{\mathfrak{Q}\delta(\mathfrak{T})\mathfrak{x},\delta(\mathfrak{T})\mathfrak{x}}\rangle=\frac{1}{2N(2N+1)}.
\end{equation}
A basic result on the steering of the canonical system (\ref{AB0})--(\ref{AB}) to zero is as follows:
\begin{theorem}\label{main}
Consider system (\ref{syst1})--(\ref{syst2}) in the canonical form (\ref{AB0})--(\ref{AB})
\begin{itemize}
	\item[A:] The matrix  $\mathfrak{Q}$ defines a common quadratic Lyapunov function for the matrices $-\mathfrak{M}$ and $\mathfrak {A+BC}.$
	\item[B:] Equation (\ref{condu}) defines ${\mathfrak T}={\mathfrak T}(\mathfrak{x})$ uniquely.
	\item[C:] Control (\ref{u}) is bounded: $|{\mathfrak u}|\leq\frac{1}{2}$.
	\item[D:] Control (\ref{u}) brings the point $\mathfrak{x}$ to $0$ in time ${\mathfrak T}(\mathfrak{x})$.
	\item[E:] The matrix $\mathfrak{Q}$ is integer and even: all elements of the matrix are even integers.
	\end{itemize}
\end{theorem}

Numerical experiments support the conjecture that strengthens the statement E of the theorem, namely, that all elements of the matrix $\mathfrak Q$ are divisible by $\mathfrak{Q}_{11}$. 
For example, for cases of one and two oscillators, the matrix has the form
\begin{equation}\label{Q2}
	\mathfrak Q_1=6\times
	\left( {\begin{array}{*{20}c}
		{1}	&	{ -2}	\\
		{-2}	&	{6}	\\
	\end{array}} \right), \quad
	\mathfrak Q_2=20\times
	\left( {\begin{array}{*{20}c}
		{1}	&	{ -9}		&	{21}		&	{-14}		\\
		{-9}	&	{111}		&	{-294}	&	{210}	\\
		{21}	&	{ -294}	&	{840}	&	{-630}	\\
		{-14}	&	{210}	&	{-630}	&	{490}	\\
	\end{array}} \right).
\end{equation}

\section{Matching of controls}\label{comb}

In Section \ref{terminal}, we solved the problem of a local feedback control that works for sufficiently small initial conditions. 
The switching to this type of control should occur at the boundary of a  domain \textit{invariant} with respect to the phase flow. 
Consider the invariant domains of the form
\begin{equation}\label{theta}
	G_\Theta=\{{\mathfrak x}: {\mathfrak T}({\mathfrak x})\leq\Theta\}=\{{\mathfrak x}: \langle\mathfrak{Q}\delta(\Theta){\mathfrak x},\delta(\Theta){\mathfrak x}\rangle\leq1\}.
\end{equation}
The invariant domain $G_\Theta$ should satisfy two conditions:
\begin{itemize}\label{conditionsG}
	\item[A:] The domain $G_\Theta$ contain the ``inefficiency'' domain $\{\rho(x)\leq UC(\underline\omega)\}$ of the preceding control,
	\item[B:] The domain  $G_\Theta$ is contained in the strip $\{|Cx|\leq1/2\}$, where $C$ is the matrix (\ref{C}).
\end{itemize}

Condition B allows us to use at the terminal stage controls ${\mathfrak u}$ which are less than 1/2 in absolute value.
Condition A reduces to the fact that set $UC(A,B)\Omega$, where $C(A,B)$ is a constant, is contained in $G_\Theta$. 
In other words, the following inequality should be fulfilled for the support functions:
\begin{equation}\label{condA}
	UC(A,B)H_\Omega(D^*p)\leq\langle{\delta(\Theta)^{-1}\mathfrak{q}\delta(\Theta)^{-1}p,p}\rangle^{1/2},
\end{equation}
where $D$ is the matrix (\ref{d}). 
It is clear that for sufficiently small $U$ the inequality holds true.

Condition B says that the value of the support function of the ellipsoid $G_\Theta$ at the vector ${D^*}^{-1}C$ does not exceed $1/2$ in absolute value. 
In other words,
\begin{equation}\label{condB}
	\left\langle\delta(\Theta)^{-1}\mathfrak{q}\delta(\Theta)^{-1}{D^*}^{-1}C,{D^*}^{-1}C\right\rangle^{1/2}\leq1/2.
\end{equation}
Of course, this inequality holds for a sufficiently small $\Theta$. 
Once $\Theta$ is chosen, we have to choose the bound $U$ for the control at the second stage in accordance with Inequality (\ref{condA}). 
Then Conditions A and B are met. 
The switching to the final stage happens upon arrival at the boundary $\{(\mathfrak{Q}\delta(\Theta){\mathfrak x},\delta(\Theta){\mathfrak x})=1\}$ of $G_\Theta$.

\section{Single linear oscillator}\label{toy}

Illustrate our preceding arguments in the simplest case of a single linear oscillator described by the equation ($\omega=1$)
\begin{equation}\label{n=12model}
\begin{array}{l}
	\dot{x}=y\\
	\dot{y}=-x+{u},
	\quad 
	|u|\leq1.
\end{array}
\end{equation}

\begin{figure}
\begin{center}
\begin{minipage}[h]{0.5\linewidth}
\center{\resizebox{1.0\columnwidth}{!}{\includegraphics{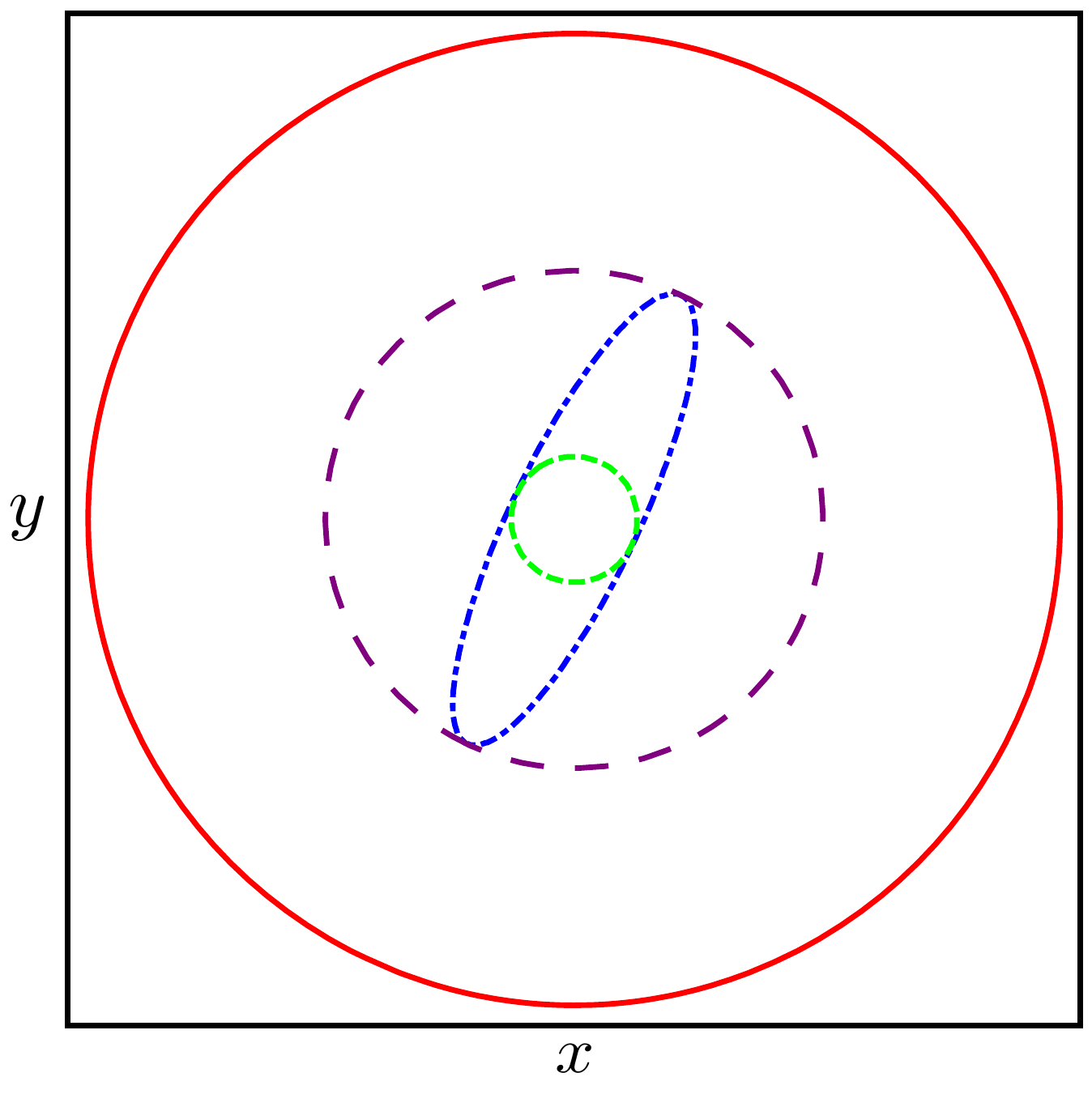}}}
\end{minipage}
\end{center}
\caption
{
	The structure of the suggested control for the case of single oscillator.
} 
\label{fig:4}
\end{figure}

We perform the separation of the phase space  $\mathbb{R}^2$ in three regions (see Fig. \ref{fig:4}). 
``Basic'' one is the outer region of the circle $\mathbb{B}_2$ with  radius $2$ (circle bounded by solid line on Fig. \ref{fig:4}), 
where the dry friction type control $u=-\sign(y)$ is applied.  
In principle, it is possible to use the $\mathbb{B}_r$ of any radius $r>1$. 
Substantially different control,
\begin{equation}
	u(x,y)=x+6\mathfrak{T}^{-2}x-3\mathfrak{T}^{-1}y
\end{equation}
is applied within a  zone close to zero (ellipse bounded by the dot-dashed line on Fig. \ref{fig:4}). 
Here $\mathfrak{T}$ is the function $(x,y)$ defined by Eq. (\ref{condB}). 
In our case, it has the form
\begin{equation}\label{conduu}
	\mathfrak{T}^{-2}6y^2-\mathfrak{T}^{-3}24xy+\mathfrak{T}^{-4}36x^2=1/6.
\end{equation}
The close to zero zone $G_\Theta$, where this control is used, 
is the interior of the ellipse $\Theta^{-2}6y^2-\Theta^{-3}24xy+\Theta^{-4}36x^2=1$ with parameter $\Theta=3^{1/4}$ found from the condition (\ref{condB}).

The ellipse contains the disk $\mathbb{B}_{\Lambda}$ (circle bounded by the dotted line on Fig. \ref{fig:4}) with radius $\Lambda=(\lambda_{\rm max})^{-1/2}=0.26253\dots$, 
where $\lambda_{\rm max}$ is the largest eigenvalue of the quadratic form (\ref{conduu}). 
If at the first stage we would use the circle $\mathbb{B}_r,\, r>1$ instead of $\mathbb{B}_2$, the parameter $U$ should be $\Lambda/r$.

\section*{Conclusion}

In our work a constructive method  has been demonstrated for the design of a feedback control for a system of arbitrary large number  of oscillators. 
In this paper, we have described computational aspects of the suggested approach to design of the control and, in particular, 
have shown that the design of the control in high energy zones reduces to solution of equation (\ref{approx3}), 
which essentially coincides with the well-studied problem of maximization a linear form on a convex hypersurface in $\mathbb{R}^N$. 
As an illustrative example, the case of a single linear oscillator is considered.

\section*{Acknowledgements}

This work is supported by Russian Foundation for Basic Research (grants 14-08-00606 and 14-01-00476).


\begin{thebibliography}{10}

\bibitem{pont}\label{pont}
{\sc L.S. Pontryagin, V.G. Boltyanskii, R.V. Gamkrelidze, and E.F. Mishchenko}.
{\em The mathematical theory of optimal processes}
(Interscience Publishers John Wiley \& Sons, Inc. New York--London, 1962).

\bibitem{kalman}\label{kalman}
{\sc R.E. Kalman}.
{\em On the general theory of control systems},
in Proceedings of the First IFAC World Congress, Moscow, {\bf 1} (1960), 481--492
(Butterworths, London, 1960).

\bibitem{ovseev6}\label{ovseev6}
{\sc  A.I. Ovseevich}. 
{\em Limit behaviour of attainable and superattainable sets}, 
{\em in Proceedings of the Conference on Modeling, Estimation and Filtering of Systems with Uncertainty} (September 3--7, Sopron, Hungary),
{\href{http://dx.doi.org/10.1007/978-1-4612-0443-5_21}{pp. 324--333, 1990}}.

\bibitem{ovseev}\label{ovseev}
{\sc E.V. Goncharova and A.I. Ovseevich}. 
{\em Comparative analysis of the asymptotic dynamics of reachable sets to linear systems}, 
Journal of Computer and Systems Sciences International, 
{\href{http://dx.doi.org/10.1134/S1064230707040016}{{\bf 46} (2007), 4, 505--513}}.

\bibitem{ovseev2}\label{ovseev2}
{\sc A.I. Ovseevich}. 
{\em Singularities of attainable sets}, 
Russian Journal of Mathematical Physics, 
{\bf 5} (1998), 3, 389--398.

\bibitem{anan}\label{anan}
{\sc I.M. Anan'evskii, N.V. Anokhin, and A.I. Ovseevich}.
{\em Bounded feedback controls for linear dynamic systems by using common Lyapunov functions}, 
Doklady Mathematics, 
{\href{http://dx.doi.org/10.1134/S106456241005039X}{{\bf 82} (2010), 2, 831--834}}.

\bibitem{anan2}\label{anan2}
{\sc I.M. Anan'evskii, N.V. Anokhin, and A.I. Ovseevich}.
{\em Design of bounded feedback controls for linear dynamical systems by using common Lyapunov functions},
Theoretical and Applied Mechanics Letters, 
{\href{http://dx.doi.org/10.1063/2.1101301}{{\bf 1} (2011), 013001}}.

\bibitem{korob}\label{korob}
{\sc  A.E. Choque Rivero, V.A. Skorik, and V.I. Korobov}.  
{\em Controllability function as the time of motion I}, 
Mathematical Physics, Analysis, Geometry (in Russian), 
{\href{http://www.mathnet.ru/php/archive.phtml?wshow=paper&jrnid=jmag&paperid=201&option_lang=eng}{{\bf 11} (2004), 2, 208--225}},
{\href{http://arxiv.org/abs/1509.05127}{arXiv:1509.05127}}.

\bibitem{chern}\label{chern}
{\sc F.L. Chernousko}. 
{\em On the construction of a bounded control in oscillatory systems}, 
Journal of Applied Mathematics and Mechanics, 
{\href{http://dx.doi.org/10.1016/0021-8928(88)90028-7}{{\bf 52} (1988), 4, 426--433}}.

\bibitem{ovseev3}\label{ovseev3}
{\sc A.I. Ovseevich}. 
{\em Complete controllability of linear dynamic systems},
Journal of Applied Mathematics and Mechanics,
{\href{http://dx.doi.org/10.1016/0021-8928(89)90119-6}{{\bf 53} (1989), 5, 665--668}}.

\bibitem{ovseev4}\label{ovseev4}
{\sc A.I. Ovseevich and A.K. Fedorov}.
{\em Feedback bounded control for a system of oscillators},
Vestnik of Lobachevsky State University of Nizhni Novgorod, 
{\href{http://www.vestnik.unn.ru/nomera?anum=6280}{{\bf 3}, 1, (2013), 278--283 [In Russian]}}.

\bibitem{ovseev5}\label{ovseev5}
{\sc A.I. Ovseevich and A.K. Fedorov}. 
{\em Asymptotically optimal feedback control for a system of linear oscillators}, 
Doklady Mathematics, 
{\href{http://dx.doi.org/10.1134/S106456241305013X}{{\bf 88} (2013), 2, 613--617}}.

\bibitem{fedorov}\label{fedorov}
{\sc A.K. Fedorov and A.I. Ovseevich}. 
{\em Asymptotic control theory for a system of linear oscillators}, 
{\href{http://arxiv.org/abs/1308.6090}{arXiv:1308.6090}}.

\bibitem{chern2}\label{chern2}
{\sc F.L. Chernousko}. 
{\em State Estimation of Dynamic Systems}
(SRC Press, Boca Raton, 1994).

\bibitem{brun}\label{brun}
{\sc P. Brunovsky}.  
{\em A classification of linear controllable systems},
{\href{http://dml.cz/bitstream/handle/10338.dmlcz/125221/Kybernetika_06-1970-3_2.pdf}{Kybernetika, {\bf 6} (1970), 173--188}}.

\end{thebibliography}
\end{document}